\documentclass{article}

\usepackage{amssymb}
\usepackage{amsmath}
\usepackage{amsthm}
\usepackage{amscd}
\usepackage{amsopn}

\theoremstyle{definition}
\newtheorem{thm}{Theorem}[section]
\newtheorem{defn}[thm]{Definition}
\newtheorem{lem}[thm]{Lemma}

\newtheorem{eg}[thm]{Example}

\newtheorem{algo}[thm]{Algorithm}

\def\K{{K}}

\def\st{\ {\rm s.t.} \ }

\def \q{~~~~~~}

\def \Zzero{{\mathbb Z}_{\geq 0}}
\def \Z{{\mathbb Z}}
\def \N{{\mathbb N}}

\def \LM{\mathrm{LM}_{<}}

\def \LE{\mathrm{LE}_{<}}
\def \EXP{\mathrm{Exps}}

\def \LEr{\mathrm{LE}_{<_s}}
\def \LMr{\mathrm{LM}_{<_r}}

\def \LTr{\mathrm{LT}_{<_r}}
\def \LEr{\mathrm{LE}_{<_r}}
\def \Restr{\mathrm{Rest}_{<_r}}
\def \Ann{\mathrm{Ann}}
\def \ord{\mathrm{ord}}

\def \inn{\mathrm{in}}

\def \nf{\mathrm{NF}}
\def \Dformal{\widehat{\mathcal{D}}}
\def \dx{\partial_x}
\def \dy{\partial_y}
\def \dz{\partial_z}
\def \d1{\partial_1}
\def \d2{\partial_2}
\def \d3{\partial_3}

\def \d{\partial}

\def \p{\partial}
\def\Proof{{\it(Proof)}~~}

\title{An Algorithm Computing the Local $b$ Function by an 
Approximate Division Algorithm in $\Dformal$}
\author{Nakayama, Hiromasa}

\begin{document}

\maketitle

\section{Introduction}

Let $\K$ be a field of characteristic $0$.
We denote the differential operator $\frac{\d}{\d x_i}$ by $\d_i$.
The ring of differential operators with polynomial coefficient 
$$D = \left\{\sum_{\beta \in \N^n} a_{k,\beta}(x) 
 \d^{\beta} ~|~ a_{k,\beta}(x) \in \K[x] \right\}$$
is denoted by $D$, and 
that with formal power series coefficient 
$$\Dformal = \left\{\sum_{\beta \in \N^n} a_{k,\beta}(x) 
 \d^{\beta} ~|~ a_{k,\beta}(x) \in \K[[x]] \right\}$$
is denoted by $\Dformal$.
Here $x = (x_1, \cdots, x_n)$ and $\d = (\d_1, \cdots, \d_n)$.

Division theorems are fundamental in several construction for $D$-modules.
Castro gave a constructive division theorem 
in $\Dformal$, 
which gives unique quotients and remainder(\cite{CastroPhD}).
However, the division procedure needs infinite reductions.
On the other hand, Granger, Oaku and Takayama gave 
a division algorithm for algebraic data in $\Dformal$ (\cite{GO},\cite{GOT}),
which is an analogous algorithm with the Mora division algorithm 
in the polynomial ring (\cite{Singular}).
We call this division the Mora division algorithm in $D$.
The division algorithm stops in finite steps, but the remainder 
is not completely reduced.
By using the division algorithm, 
we can get a Gr\"{o}bner basis and solve the ideal membership 
problem for algebraic data in $\Dformal$.

The local $b$ function of a polynomial $f \in \K[x]$ is 
the minimum degree monic polynomial $b(s) \in \K[s]$ which satisfies  
$$ \exists P \in \Dformal[s] \text{ such that } Pf^{s+1} = b(s) f^s.$$
Oaku gave algorithms computing the local $b$ function of a given polynomial
by using a Gr\"{o}bner basis method in $D$ 
(\cite{Oakubfunc1},\cite{Oakubfunc2},\cite{OTbfunc}).
In this paper,
we propose a new algorithm to compute the local $b$ function.
We use the Mora division algorithm in $D$    
and an approximate division algorithm in $\Dformal$, which gives
an approximation of the remainder by Castro's 
division in $\Dformal$.

\section{Division Theorem and Approximate Division Algorithm in the Ring of 
Power Series $\K[[x]]$}

There are several kinds of division theorem in $\K[[x]]$.
Among them, we need a division theorem which gives unique quotients and 
unique remainder for our approximate division algorithm in $\Dformal$.
We start with explaining this kind of division theorem.
As to details about the division theorem including a history, we refer
\cite{CG}.

We define a monomial order $<_r$ on monomials in $\K[[x]]$ by the following $1 \times n$
matrix and a term order $<$ as the tie-breaker.
\begin{center}
 \begin{tabular}{ccc} 
  $x_1$ & $\cdots$ & $x_n$ \\
  \hline
  $-1$  & $\cdots$ & $-1$  \\
  $($   & term order $<$ & $)$
 \end{tabular} 
\end{center}
In other words, we define the order $<_r$ as 
$$
x^{\alpha} <_r x^{\beta} \Leftrightarrow
\begin{cases}
-\alpha_1- \cdots -\alpha_n < -\beta_1- \cdots -\beta_n \text{ or } \\
(-\alpha_1- \cdots -\alpha_n = -\beta_1- \cdots -\beta_n \text{ and } x^{\alpha} < x^{\beta})
\end{cases}
$$
Since the order is not a well order, ordinary division algorithm 
does not generally stop in finite steps.

Let $\alpha(1), \cdots, \alpha(s) \in (\Zzero)^n$.
For $(\alpha(1), \cdots, \alpha(s))$, we define a partition 
$\Delta^1, \cdots, \Delta^s, \overline{\Delta}$ of $(\Zzero)^n$ as  
\begin{align*}
\Delta^1 &= \alpha(1) + (\Zzero)^n, \\
\Delta^2 &= (\alpha(2) + (\Zzero)^n) \setminus \Delta^1, \cdots \\
\Delta^s &= (\alpha(s) + (\Zzero)^n) \setminus (\Delta^1 \cup \cdots \cup \Delta^{s-1}), \\
\overline{\Delta} &= (\Zzero)^n \setminus (\Delta^1 \cup \cdots \cup \Delta^s).
\end{align*}

Let $f \in \K[[x]]$. 
We define the following notations.  
\begin{align*}
\LMr(f) :& \text{ leading monomial of } f \text{ with respect to } <_r \\
\LTr(f) :& \text{ leading term of } f \text{ with respect to } <_r \\
\LEr(f) :& \text{ leading exponent of } f \text{ with respect to } <_r \\
\EXP(f) :& \text{ set of exponents appearing in } f \\
\Restr(f) :& \text{ } f - \LTr(f)
\end{align*}
For example, let $f = 3 x_1 + x_1 x_2$, then 
$\LMr(f) = x_1, \LTr(f) = 3x_1, \LEr(f) = (1,0), 
\EXP(f) = \left\{(1,0), (1,1)\right\}$, 
and $\Restr(f) = x_1 x_2$.

\begin{lem}
(Division by monomials in $\K[[x]]$)
\label{monodiv}

Let $\alpha(1), \cdots, \alpha(s) \in (\Zzero)^n$ 
and $f \in \K[[x]]$ and $\Delta^1, \cdots, \Delta^s, \overline{\Delta}$ be 
the partition of $(\Zzero)^n$ with respect to $(\alpha(1), \cdots, \alpha(s))$.
Then,
there exist quotients $q_1, \cdots, q_s \in \K[[x]]$ and remainder 
$r \in \K[[x]]$ which satisfies the following conditions :
\begin{align*}
&  f = q_1 x^{\alpha(1)} + \cdots + q_s x^{\alpha(s)} + r, \\
&  \alpha(i) + \EXP(q_i) \subset \Delta^i, \\
&  \EXP(r) \subset \overline{\Delta}.
\end{align*}
Especially, $(q_1, \cdots, q_s, r)$ is uniquely determined by $(f, x^{\alpha(1)}, \cdots, x^{\alpha(s)}, <_r)$.
\end{lem}

\begin{thm}
(Division theorem in $\K[[x]]$)
\label{power-series-division}

Let $f \in \K[[x]]$,  
$g_1, \cdots, g_s \in \K[[x]]$, and
$\Delta^1, \cdots, \Delta^s, \overline{\Delta}$ be the partition of $(\Zzero)^n$ with respect to 
$(\LEr(g_1), \cdots, \LEr(g_s))$.
Then, 
there exist quotients $q_1, \cdots, q_s \in \K[[x]]$ and 
remainder $r \in \K[[x]]$ which satisfies the following conditions.
\begin{align*}
& f = q_1 g_1 + \cdots + q_s g_s + r \\
& \LEr(g_i) + \EXP(q_i) \subset \Delta^i \\
& \EXP(r) \subset \overline{\Delta} 
\end{align*}
Especially, $(q_1, \cdots, q_s, r)$ is uniquely determined by $(f, g_1, \cdots, g_s, <_r)$.
\end{thm}

We show a procedure to obtain $q_1, \cdots, q_s, r$.

\begin{center}
\begin{tabular}{l}
$K[[x]]$-division$(f, g_1, \cdots, g_s)$\\
\hline
Input : $f \in \K[[x]], g_1, \cdots, g_s \in \K[[x]]$\\
Output : $q_1, \cdots, q_s, r \in \K[[x]]$ which satisfies \\
~~~~~~~~~ $f = q_1 g_1 + \cdots q_s g_s + r$\\
~~~~~~~~~ $\LEr(g_i) + \EXP(q_i) \subset \Delta^i$\\
~~~~~~~~~ $\EXP(r) \subset \overline{\Delta}$ \\ 
\hline
$1:$\q $F \leftarrow f, q_i \leftarrow 0 ~(1 \leq i \leq s), r \leftarrow 0$\\
$2:$\q while ($F \neq 0$) \{\\
$3:$\q \q $[q',r'] \leftarrow \text{mono-div}(F, g_1, \cdots, g_s)$ \\
$4:$\q \q $q_i \leftarrow q_i + q_i' ~(1 \leq i \leq s)$ \\
$5:$\q \q $r \leftarrow r + r'$ \\
$6:$\q \q $F \leftarrow -\sum_{i=1}^s q_i' \Restr(g_i)$\\
$7:$\q \}\\
$8:$\q return [$q,r$] \\
\hline
\end{tabular}
\end{center}

Here, the procedure mono-div$(F, G)$ at the 3rd line returns 
$q_1', \cdots, q_s', r' \in \K[[x]]$ which satisfies 
\begin{align*}
& F = q_1' \LTr(g_1) + \cdots + q_s' \LTr(g_s) + r' \\
& \LEr(g_i) + \EXP(q_i') \subset \Delta^i \\
& \EXP(r') \subset \overline{\Delta} 
\end{align*}
Namely, divide $F$ by $\LTr(G)$ (Lemma \ref{monodiv}).

This procedure generally does not stop in finite steps.
Therefore this procedure is not an algorithm in a strict sence.
If input $f$ and $g_i$ are polynomials, 
we can give an algorithm which returns approximate quotients and approximate
remainder which are correct up to given total degree $N$.
Approximate algorithm of this kind has been discussed by several authors, 
for example, see \cite{GBPowerSeries}.
After the manner of the approximate algorithm, 
we will give an algorithm which gives approximations of quotients and remainder 
by $K[[x]]$-division. 

\begin{algo}
($K[[x]]$-approximate-division)
\label{power-series-app-div} 
\begin{center}
\begin{tabular}{l}
$K[[x]]$-approximate-division$(f, g_1, \cdots, g_s)$\\
\hline
Input : $f \in \K[x], g_1, \cdots, g_s \in \K[x], N \in \Z_{>0}$ \\
Output : $\overline{q_1}, \cdots, \overline{q_s}, \overline{r}, \overline{F} \in \K[x]$ which satisfies the following
conditions. \\ 
~~~~~~~~~~ $f = \overline{q_1} g_1 + \cdots + \overline{q_s} g_s + \overline{r} + \overline{F}$\\
~~~~~~~~~~ $\overline{F} \neq 0 \Rightarrow |\LEr(\overline{F})| \geq N$\\ 
~~~~~~~~~~ Terms of $\overline{q_i}$ whose total degree is smaller than $N - |\LEr(g_i)|$ agree \\ 
~~~~~~~~~~ with those of $q_i$.\\
~~~~~~~~~~ Terms of $\overline{r}$ whose total degree is smaller than $N$ agree with those of $r$.\\
Here $|\alpha|$ is $\sum_{i=1}^n \alpha_i$, where $\alpha = (\alpha_1, \cdots, \alpha_n)$, \\ 
and $q_i, r$ are quotients and remainder of $K[[x]]$-division($f,G$).\\
\hline
$1:$ \q $\overline{F} \leftarrow f, \overline{q_i} \leftarrow 0 ~(1 \leq i \leq s), \overline{r} \leftarrow 0$\\
$2:$ \q while ($\overline{F} \neq 0$ and $|\LEr(\overline{F})| < N$) \{\\
$3:$ \q \q $[q',r'] \leftarrow \text{mono-div}(\overline{F}, G)$ \\
$4:$ \q \q $\overline{q_i} \leftarrow \overline{q_i} + q_i' ~(1 \leq i \leq s)$ \\
$5:$ \q \q $\overline{r} \leftarrow \overline{r} + r'$ \\
$6:$ \q \q $\overline{F} \leftarrow -\sum_{i=1}^s q_i' \Restr(g_i)$\\
$7:$ \q \}\\
$8:$ \q return [$\overline{q}, \overline{r}, \overline{F}$] \\
\hline
\end{tabular}
\end{center}
\end{algo}

\Proof
We denote variables $\overline{F}, \overline{g_i}, \overline{r}, q_i', r'$ at the end of the execution of the $k$-th while-loop
by $\overline{F}^{(k)}, \overline{g_i}^{(k)}, \overline{r}^{(k)}, q_i'^{(k)}, r'^{(k)}$.
From the result of mono-div$(\overline{F}^{(k-1)}, G)$, it holds that
\begin{align}
&\overline{F}^{(k-1)} = \sum_{i=1}^s q_i'^{(k)} g_i + r'^{(k)} +
 \overline{F}^{(k)} \label{KtoKplusOneE}\\
&\LMr(q_i'^{(k)} g_i) \leq_r \LMr(\overline{F}^{(k-1)}), 
\LMr(r'^{(k)}) \leq_r \LMr(\overline{F}^{(k-1)})
\label{ineq-monodiv}
\end{align}

We will prove that the algorithm stops in finite steps.
From (\ref{KtoKplusOneE}) and (\ref{ineq-monodiv}), we get
$$ \LMr(\overline{F}^{(k)}) \geq_r \LMr(\overline{F}^{(k+1)}) \geq_r \LMr(\overline{F}^{(k+2)}) \geq_r \cdots .$$
From the definition of the monomial order $<_r$, we get
$$ |\LEr(\overline{F}^{(k)})| \leq |\LEr(\overline{F}^{(k+1)})| \leq |\LEr(\overline{F}^{(k+2)})| \leq \cdots .$$
Since the number of monomials which has the same total degree is finite, 
there exists $k$ such that $|\LEr(\overline{F}^{(k)})| \geq N$.
This is the condition which stops the while loop.
Therefore the algorithm stops in finite steps.

Next, we will prove that the returned $\overline{q_i}, \overline{r}$ satisfy the conditions of output. 
At the end of the execution of the $k$-th while-loop, it holds that  
$$
f = \sum_{i=0}^s \overline{q_i}^{(k)} g_i + \overline{r}^{(k)} + \overline{F}^{(k)} 
$$
Since $q'^{(k)}, r'^{(k)}$ is the result of mono-div$(\overline{F}^{(k-1)}, G)$, it holds that 
\begin{align*}
&\LEr(g_i) + \EXP(q_i'^{(k)}) \subset \Delta^i\\
&\EXP(r'^{(k)}) \subset \overline{\Delta}
\end{align*}
From these properties, we have
$$\LMr(\overline{F}^{(k-1)}) = \max(\LMr(q_1'^{(k)}g_1), \cdots, \LMr(q_s'^{(k)}g_s), \LMr(r'^{(k)})) $$
That is 
$$ 
\LMr(q_i'^{(k)} g_i) \leq_r \LMr(\overline{F}^{(k-1)}), \LMr(r'^{(k)}) \leq_r \LMr(\overline{F}^{(k-1)}) $$
So we get 
\begin{align*}
&|\LEr(q_i'^{(k)})| \geq |\LEr(\overline{F}^{(k-1)})| - |\LEr(g_i)|,\\ 
&|\LEr(r'^{(k)})| \geq |\LEr(\overline{F}^{(k-1)})|.
\end{align*}
Since $\overline{q_i}^{(k)} = \overline{q_i}^{(k-1)} + q_i'^{(k)}$ and $\overline{r}^{(k)} = \overline{r}^{(k-1)} + r'^{(k)}$,
terms of $\overline{q_i}^{(k)}$ whose total degree is smaller than 
$|\LEr(\overline{F}^{(k-1)})| - |\LEr(g_i)|$ does not change, 
and terms of $\overline{r}^{(k)}$ whose total degree is smaller than $|\LEr(\overline{F}^{(k-1)})|$
does not change.
Therefore $\overline{q_i}, \overline{r}$ satisfy the condition.
$\blacksquare$

\begin{eg}
(Example of $K[[x]]$-approximate-division)

The example of the case of $f = x, G={1+x}, N=5$ (execute $K[[x]]$-approximate-division$(x, 1+x, 5)$) 
$$
\begin{matrix}
\text{mono-div}         & q' & r'& \overline{q}          & \overline{r} & \overline{F}         \\
\hline
                        &    &   & 0                & 0       & x               \\
x    = x \cdot 1 + 0    & x  & 0 & x                & 0       & -x \cdot x      \\
-x^2 = -x^2 \cdot 1 + 0 &-x^2& 0 & x-x^2            & 0       & -(-x^2) \cdot x \\
x^3  = x^3 \cdot 1 + 0  & x^3& 0 & x-x^2+x^3        & 0       & -x^3 \cdot x    \\
-x^4 = -x^4 \cdot 1 + 0 &-x^4& 0 & x-x^2+x^3-x^4    & 0       & -(-x^4) \cdot x \\
x^5  = x^5 \cdot 1 + 0  & x^5& 0 & x-x^2+x^3-x^4+x^5& 0       & -x^5 \cdot x 
\end{matrix}
$$
The output of the approximate division algorithm is 
$$ x = (x-x^2+x^3-x^4+x^5) \cdot (1+x) - x^6.$$
In this case, the quotients and the remainder by $K[[x]]$-division are
\begin{align*}
 x  &= \frac{x}{1+x} \cdot (1+x) + 0 \\
    &= (x-x^2+x^3-x^4+x^5- \cdots) \cdot (1+x) + 0 
\end{align*}
\end{eg}

\section{Division Theorem and Approximate Division Algorithm in $\widehat{\mathcal{D}}[y]$}

Castro gave a constructive division procedure in 
$\Dformal$~(\cite{CastroPhD}).
His division procedure returns unique quotients and remainder.
In this section, we give an algorithm which gives approximations of 
the quotients and the remainder by his division.
The approximate quotients and remainder 
are correct up to a given total degree.  

Let $y$ be a parameter, and
$\xi$ be a commutative variable which stands for $\p$.
We define a monomial order $<$ on $\Dformal[y]$ 
by using the following $2 \times (2n + 1) $ matrix and the tie-breaker $<_1$. 
\begin{center}
 \begin{tabular}{ccccccc}
 $x_1$ & $\cdots$ & $x_n$ & $y$ & $\xi_1$ & $\cdots$ & $\xi_n$ \\
 \hline
 $0$   & $\cdots$ &   $0$ & $1$ &     $1$ & $\cdots$ &     $1$ \\
 $-1$  & $\cdots$ &   $-1$& $0$ &     $0$ & $\cdots$ &     $0$ \\
 $($   & \multicolumn{5}{c}{term order ~ $<_1$}       & $)$
 \end{tabular}
\end{center}
And we define an other monomial order $<_r$ on $\Dformal[y]$ 
by using the following $1 \times (2n + 1) $ matrix and tie-breaker $<_1$.
\begin{center}
 \begin{tabular}{ccccccc}
 $x_1$ & $\cdots$ & $x_n$ & $y$ & $\xi_1$ & $\cdots$ & $\xi_n$ \\
 \hline
 $-1$   & $\cdots$ & $-1$ & $-1$&    $-1$ & $\cdots$ &    $-1$ \\
 $($   & \multicolumn{5}{c}{term order $<_1$}       & $)$
 \end{tabular}
\end{center}
We define a weight vector $e$ for monomials $\left\{x^{\alpha} \xi^{\beta} y^{\gamma}\right\}$ as
\begin{center}
 \begin{tabular}{ccccccccc}
    & $x_1$ & $\cdots$ & $x_n$ & $y$ & $\xi_1$ & $\cdots$ & $\xi_n$ \\
$($ & $0$   & $\cdots$ & $0$   & $1$&    $1$ & $\cdots$ &    $1$ & $)$
 \end{tabular}
\end{center}

If $|\beta| + |\gamma| = |\beta'| + |\gamma'|$ 
(in other words, monomials have the same order with respect to the weight vector $e$), 
then it holds that
$$x^{\alpha} \xi^{\beta} y^{\gamma} < x^{\alpha'} \xi^{\beta'} y^{\gamma'}
 \Leftrightarrow
 x^{\alpha} \xi^{\beta} y^{\gamma} <_r x^{\alpha'} \xi^{\beta'} y^{\gamma'}$$
Therefore, we get
$$ \LM(P) = \LMr(\inn_e(P))$$
Here $\inn_e(P)$ is the initial part of $P$ with respect to the weight vector $e$.
For example, when $P = x_1 \p_1^2 + x_1^2 \p_1^2 + x_1 \p_1 + 1$, 
we have $\inn_e(P) = x_1 \xi_1^2 + x_1^2 \xi_1^2$. 
And $\LMr(\inn_e(P)) = x_1 \xi_1^2, \LM(P) = x_1 \xi_1^2$ 

We review Castro's division procedure in $\Dformal[y]$.
Of course, 
since the monomial order $<$ is not a well order, 
the procedure needs infinite reductions.

\begin{thm}
(Division theorem in $\Dformal[y]$, \cite{CastroPhD})
\label{D-division}

Let $P, P_1, \cdots, P_s \in \Dformal[y]$, and $\Delta^1, \cdots,
 \Delta^s, \overline{\Delta}$ be the partition of $\Zzero^{2n+1}$ with respect
 to $(\LE(P_1), \cdots, \LE(P_s))$. 
There exist $Q_1, \cdots , Q_s, R \in \Dformal[y]$ which satisfies the following conditions.
\begin{align*}
&P = Q_1 P_1 + \cdots + Q_s P_s + R \\
&\LE(P_i) + \EXP(Q_i) \subset \Delta^i \\
&\EXP(R) \subset \overline{\Delta} 
\end{align*}
Especially, $(Q_1, \cdots, Q_s, R)$ are uniquely determined by $(P, P_1,
\cdots, P_s, <)$.
\end{thm}

We show procedure for the division theorem.

\begin{center}
\begin{tabular}{l}
$\Dformal$-division($P, P_1, \cdots, P_s$)\\
\hline
Input : $P, P_1, \cdots, P_s \in \Dformal[y]$ \\
Output : $Q_1, \cdots, Q_s, R \in \Dformal[y]$ which satisfies \\ 
~~~~~~~~~~~~ $P = Q_1 P_1 + \cdots + Q_s P_s + R$ \\
~~~~~~~~~~~~ $\LE(P_i) + \EXP(Q_i) \subset \Delta^i$\\
~~~~~~~~~~~~ $\EXP(R) \subset \overline{\Delta}$ \\
\hline
$1:$~~ \q $Q_i \leftarrow 0 ~~ (1 \leq i \leq s), R \leftarrow P$ \\
$2:$~~ \q $m_0 \leftarrow \ord_e(R)$ \\
$3:$~~ \q for ($k \leftarrow m_0; k \geq 0; k \leftarrow k - 1$) \{ \\
$4:$~~ \q \q $r \leftarrow$ (total symbol of the part of $R$ whose $e$-order is $k$) \\
$5:$~~ \q \q [$q_i', r'$] $\leftarrow$ $K[[x]]$-division $(r, \inn_e(P_1), \cdots, 
           \inn_e(P_s), <_r$) \\
$6:$~~ \q \q $Q_i' \leftarrow$ (replace $\xi$ with $\partial$ in $q_i'$) ~ ($1 \leq i \leq s$) \\
$7:$~~ \q \q$R \leftarrow R - \sum_{i=1}^{s} Q_i' P_i$ \\
$8:$~~ \q \q$Q_i \leftarrow Q_i + Q_i'$ ~ ($1 \leq i \leq s$)\\
$9:$~~ \q \} \\
$10:$ \q return [$Q_i, R$] \\
\hline
\end{tabular}
\end{center}

Here, $\ord_e(R)$ is $e$-order of $R$, defined by   
$$\ord_e(R) = \max\left\{e \cdot \alpha ~|~ \alpha \in \EXP(R) \right\},$$
and the total symbol of $P(x, y, \p) = \sum_{\alpha} a_{\alpha}(x,y) \p^\alpha$ 
($a_{\alpha}(x,y) \in \K[[x]][y]$) is
$P(x, y, \xi) = \sum_{\alpha} a_{\alpha}(x,y) \xi^\alpha$.
For example, when $P = (x_1 + x_1 x_2)\p_1^2 + x_1 \p_1 + 1$, 
$\ord_e(P) = 2$ and
the total symbol $P(x,\xi) = (x_1 + x_1 x_2) \xi_1^2 + x_1 \xi_1 + 1$.  

Since the procedure uses $K[[x]]$-division at the $5$th line, it does not stop in finite steps.
We suppose that inputs are algebraic data, in other words, inputs are elements in $D$.
By replacing the $K[[x]]$-division with $K[[x]]$-approximate-division 
(Algorithm \ref{power-series-app-div}), 
we can get approximations of the quotients and the remainder by
$\Dformal$-division.

\begin{algo}
\label{ddivalgo}
($\Dformal$-approximate-division)
\begin{center}
\begin{tabular}{l}
$\Dformal$-approximate-division($P, P_1, \cdots, P_s, N$)\\
\hline
Input : $P, P_1, \cdots, P_s \in D[y], N \in \Z_{>0}$\\
Output : $\overline{Q_1}, \cdots, \overline{Q_s}, \overline{R} \in D[y]$ satisfies \\
$P = \overline{Q_1} P_1 + \cdots + \overline{Q_s} P_s + \overline{R}$ \\
Terms of $\overline{Q_i}$ whose total degree is smaller than $N - |\LE(P_i)|$ agree with those of $Q_i$ \\
Terms of $\overline{R}$ whose total degree is smaller than $N$ agree with those of $R$ \\
Here, $Q_i$ and $R$ are the result of $\Dformal$-division$(P, P_1,
 \cdots, P_s)$. \\
\hline
\end{tabular}
\begin{tabular}{l}
$1:$~~ \q $ \overline{Q_i} \leftarrow 0 ~~ (1 \leq i \leq s), 
          \overline{R} \leftarrow P$ \\
$2:$~~ \q $m_0 \leftarrow \ord_e(\overline{R})$\\
$3:$~~ \q for ($k \leftarrow 0; k \leq m_0; k \leftarrow k + 1$) \\
$4:$~~ \q \q $M_k \leftarrow \max(|\LE(P_i)| + 2(k - \ord_e(P_i)) ~|~ 
 i \text{ such that } k \geq \ord_e(P_i))$ \\
$5:$~~ \q Bound $\leftarrow N + \sum_{i=0}^{m_0} M_i$\\
$6:$~~ \q for ($k \leftarrow m_0; k \geq 0; k \leftarrow k - 1$) \{\\
$7:$~~ \q \q $\overline{r} \leftarrow$ (total symbol of the part of $\overline{R}$ whose $e$-order is k and whose \\
\q \q \q \q ~ total degree is smaller than Bound)\\
$8:$~~  \q \q [$\overline{q_i'}, \overline{r'}$] $\leftarrow$ $K[[x]]$-approximate-division($
             \overline{r}, \inn_e(P_1), \cdots, \inn_e(P_s), <_r,$ \\ 
\q \q \q \q ~~~~~ Bound) \\
$9:$~~  \q \q $\overline{Q_i'} \leftarrow$ (replace $\xi$ with $\partial$ in $\overline{q_i'}$) ~ ($1 \leq i \leq s$) \\
$10:$  \q \q $\overline{R} \leftarrow \overline{R} - \sum_{i = 0}^s  \overline{Q_i'} P_i$ \\
$11:$  \q \q $\overline{Q_i} \leftarrow \overline{Q_i} + \overline{Q_i'}$ ~ ($1 \leq i \leq s$)\\
$12:$  \q \q Bound $\leftarrow$ Bound $ -M_k$\\
$13:$  \q \} \\
$14:$  \q return [$\overline{Q_i}, \overline{R}$]\\
\hline
\end{tabular}
\end{center}
\end{algo}

Since $K[[x]]$-approximate-division stops in finite steps, 
the algorithm stops in finite steps.
In the remainder of this section, we will prove the correctness of the algorithm.

We put
\begin{align*}
& D(m) = \left\{P \in D[y] ~(\text{or}~ \Dformal[y]) ~|~ \ord_e(P) = m \right\},\\
& T(i) = \left\{P \in D[y] ~(\text{or}~ \Dformal[y], \K[x, \xi], 
  \K[[x]][y, \xi]) ~|~ 
  (\text{the total degree of every term of } P) \geq  i \right\}. 
\end{align*}
It follows from the definition that  
$P - Q \in T(i)$  is equivalent to that $P$ and $Q$ are the same up to total
degree $i - 1$.
From the Leibnitz rule, we get the following property of the multiplication
of approximate elements. 

\begin{lem}
\label{app-mul-deg}
(Multiplication of approximate elements)
$$(T(i) \cap D(m)) \cdot T(j) \subset T(i + j - 2m)$$
Here, dot $\cdot$ means multiplication in the ring of differential operators.
\end{lem}
\Proof
We suppose that $P \in T(i) \cap D(m)$ and $Q \in T(j)$.
Then we have 
$$ P = \sum_{|\beta| \leq m, |\alpha| + |\beta| \geq i} a_{\alpha \beta} 
   x^\alpha \partial^\beta ~,~ 
   Q = \sum_{|\gamma| + |\delta| \geq j} b_{\gamma \delta} 
   x^\gamma \partial^\delta ~~(a_{\alpha \beta}, b_{\gamma \delta} \in \K) $$
To compute the multiplication $PQ$, we use the Leibnitz rule.
We can get the following total symbol for $PQ$.
\begin{align*}
(PQ)(x, \xi) &= \sum_\nu \frac{1}{\nu!} 
                \frac{\partial^{|\nu|}}{\partial \xi^\nu} 
                (\sum_{(\alpha, \beta) \in \EXP(P)} a_{\alpha \beta} x^\alpha
                 \xi^\beta)
		\frac{\partial^{|\nu|}}{\partial x^\nu}
		(\sum_{(\gamma, \delta) \in \EXP(Q)} b_{\gamma \delta} x^\gamma
                 \xi^\delta) \\
             &= \sum_{(\alpha, \beta) \in \EXP(P), (\gamma, \beta) \in \EXP(Q)
                ,\nu \leq \beta, \nu \leq \gamma} 
                c_{\alpha \beta \gamma \delta \nu} x^{\alpha + \gamma - \nu}
                \xi^{\beta + \delta - \nu} ~~(c_{\alpha \beta \gamma \delta 
                \nu} \in \K).
\end{align*}
Here, $\nu \leq \beta$ means $\nu_i \leq \beta_i$ for all $i$.

We consider the minimum total degree of terms of $PQ$.
\begin{align}
\label{min-td}
\min(|\alpha| + |\gamma| - |\nu| + |\beta| + |\delta| - |\nu| ~|~ 
      (\alpha, \beta) \in \EXP(P), (\gamma, \delta) \in \EXP(Q), 
      \nu \leq \beta, \nu \leq \gamma) 
\end{align}
Since 
$|\alpha| + |\beta| \geq i, |\gamma| + |\delta| \geq j$ and 
$|\nu| \leq |\beta| \leq \ord_e(P) = m$ hold, 
the minimum total degree (\ref{min-td}) is more than or equal to $i + j - 2m$.
Therefore, we conclude $PQ \in T(i + j - 2m)$. 
$\blacksquare$

\begin{lem}
\label{power-series-div-deg}

Let $q_i$ and $r$ be the quotients and remainder of \\ 
$K[[x]]$-division($f, \left\{g_1, \cdots, g_s\right\}, <_r$).
If $f \in T(N)$, then $q_i \in T(N - |\LEr(g_i)|), r \in T(N)$.
\end{lem}
\Proof
It holds that
\begin{align*}
&f=q_1 g_1 + \cdots + q_s g_s + r\\
&\LEr(g_i) + \EXP(q_i) \subset \Delta^i \\
&\EXP(r) \subset \overline{\Delta} 
\end{align*}
, where $\Delta^1, \cdots, \Delta^s, \overline{\Delta}$ are the partition of $\Zzero^{2n+1}$
with respect to $(\LEr(g_1), \cdots, \LEr(g_s))$.
From these properties, we get 
\begin{align*}
\LEr(q_i g_i) \leq_r \LEr(f) \\
\LEr(r) \leq_r \LEr(f)
\end{align*}
From the definition of the monomial order $<_r$, we have  
$$ |\LEr(q_i)| \geq N - |\LEr(g_i)| \Leftrightarrow 
   q_i \in T(N - |\LEr(g_i)|) $$
and 
$$ |\LEr(r)| \geq |\LEr(f)| = N \Leftrightarrow r \in T(N).$$
$\blacksquare$

\begin{lem}
\label{power-series-div-appdiv}

Let $f, g_1, \cdots, g_s \in \K[[x]]$ and 
$\overline{f}$ be the part of $f$ whose total degree is less than $N$,
that is, $f - \overline{f} \in T(N)$.
Let $q_i$ and $r$ be the quotients and remainder of 
$K[[x]]$-division($f, g_1, \cdots, g_s, <_r$).
And let $\overline{q_i}$ and $\overline{r}$ be the quotients and remainder of   
$K[[x]]$-approximate-division($\overline{f}, g_1, \cdots, g_s, <_r, N$).
Then $q_i - \overline{q_i} \in T(N - |\LEr(g_i)|), r - \overline{r} \in T(N)$.
\end{lem}
\Proof 
Let $q_i'$ and $r'$ be the quotients and remainder of \\
$K[[x]]$-division($\overline{f}, g_1, \cdots, g_s, <_r$).
Since $K[[x]]$-division uniquely gives the quotient and remainder,
the quotients and remainder of 
$K[[x]]$-division($f-\overline{f}, g_1, \cdots, g_s, <_r$) 
are $q_i - q_i'$ and $r - r'$.
From $f-\overline{f} \in T(N)$ and Lemma \ref{power-series-div-deg}, 
$q_i - q_i' \in T(N - |\LEr(g_i)|)$ and $r - r' \in T(N)$ holds.

From the property of $K[[x]]$-approximate-division,  
$q_i' - \overline{q_i} \in T(N - |\LEr(g_i)|)$ and $r' - \overline{r} \in T(N)$ hold.
So we get $q_i - \overline{q_i} \in T(N - |\LEr(g_i)|)$ and $r - \overline{r} \in T(N)$.
$\blacksquare$ 

{\it(Proof of Algorithm \ref{ddivalgo})}
We will compare steps in the for-loop in $\Dformal$-division with 
those in $\Dformal$-approximate division.
We put
$$M_k = \max(|\LE(P_i)| + 2(k - \ord_e(P_i)) 
 ~|~ i \text{ such that } k \geq \ord_e(P_i)).$$
We suppose that  $R - \overline{R} \in T(N)$ in the $(k-1)$-th step.
Let us compare $k$-th steps in the for-loop in $\Dformal$-division and $\Dformal$-approximate-division.

\begin{center}
\begin{tabular}{l}
$\Dformal$-division\\
\hline
$r \leftarrow ($total symbol of the part of $R$ whose $e$-order is $k)$\\
$[q_i', r'] \leftarrow$ $K[[x]]$-division($r, \left\{\inn_e(P_1), \cdots, \inn_e(P_s)\right\}, 
  <_r$)\\
$Q_i' \leftarrow q_i'(x, \partial)$\\
$R_{new} \leftarrow R - \sum Q_i' P_i$ \\
\hline
\\
$\Dformal$-approximate-division\\
\hline
$\overline{r} \leftarrow ($total symbol of the part of $\overline{R}$ whose $e$-order is $k$ and total degree  \\
~~~~~ is less than $N)$\\
$[\overline{q_i'}, \overline{r'}] \leftarrow$ $K[[x]]$-approximate-division(
$\overline{r}, \left\{\inn_e(P_1), \cdots, \inn_e(P_s)\right\}, <_r, N$)\\
$\overline{Q_i'} \leftarrow \overline{q_i'}(x, \partial)$\\
$\overline{R_{new}} \leftarrow \overline{R} - \sum \overline{Q_i'} P_i$ \\
\hline
\end{tabular}
\end{center}

We will prove that $R_{new} - \overline{R_{new}} \in T(N - M_k)$ holds after these computations have been done.

It follows from Lemma \ref{power-series-div-appdiv} that 
\begin{align*}
&q_i' - \overline{q_i'} \in T(N - |\LEr(\inn_e(P_i))|) = T(N -
 |\LE(P_i)|), \\
&r' - \overline{r'} \in T(N)
\end{align*}
Remainders $r$ and $\overline{r}$ are $e$-homogeneous elements of order $k$, 
and $\inn_e(P_i)$ is $e$-homogeneous element of order $\ord_e(P_i)$.
So, if $q_i', \overline{q_i'} \neq 0$, then $q_i', \overline{q_i'}$ are $e$-homogeneous
elements of order $k - \ord_e(P_i) (\geq 0)$.
And, if $r', \overline{r'} \neq 0$, then $r', \overline{r'}$ are $e$-homogeneous 
elements of order $k$.
Therefore, if $Q_i', \overline{Q_i'} \neq 0$, then $Q_i', \overline{Q_i'}$ have the same
$e$-order $k - \ord_e(P_i)$ and are the same up to total degree 
$N - |\LE(P_i)| - 1$.
In other words, it holds that 
$$Q_i' - \overline{Q_i'} \in T(N - |\LE(P_i)|) \cap D(k - \ord_e(P_i)).$$

From Lemma \ref{app-mul-deg}, we have 
$$(T(N - |\LE(P_i)|) \cap D(k - \ord_e(P_i)))
\cdot T(0) \subset T(N - |\LE(P_i)| - 2(k - \ord_e(P_i))).$$ 
Therefore it holds that
$$Q_i' P_i - \overline{Q_i'} P_i \in 
T(N - |\LE(P_i)| - 2(k - \ord_e(P_i)))$$
Since $M_k = \max(|\LE(P_i)| + 2(k - \ord_e(P_i)) ~|~ i \text{ such that } k \geq \ord_e(P_i))$, we get 
$$\sum_{i=1}^s Q_i' P_i - \sum_{i=1}^s \overline{Q_i'} P_i \in T(N - M_k).$$
Therefore  $R_{new} - \overline{R_{new}} \in T(N - M_k)$.

The accuracy of approximation decreases by $M_k$ at each step of the for loop.
To keep the accuracy, we beforehand add $\sum_{i=0}^{m_0} M_i$ to $N$.
 
$\blacksquare$

\begin{eg}
(Example of $\Dformal$-division and $\Dformal$-approximate-division)

Let $n = 1$.
We put $P = \partial^2, P_1 = (1+x) \partial + x$.
We compare $\Dformal$-division($P, P_1$) and 
$\Dformal$-approximate-division($P, P_1, 5$).

At first, we show the procedure of $\Dformal$-division(
$P, P_1$).

\begin{center}
\begin{tabular}{l}
\hline
$R \leftarrow P, m_0 \leftarrow \ord_e(R) = 2$ \\
\hline
$r \leftarrow ($part of $R$ whose $e$-order is $2) = \xi^2$  \\
$[q_1', r'] \leftarrow$ $K[[x]]$-division($r, \inn_e(P_1), <_r$) \\
($q_1' = \frac{1}{1+x} \xi, r' = 0$) \\
$Q_1' \leftarrow \frac{1}{1+x} \p$ \\
$R \leftarrow R - Q_1' P_1 = -\p - \frac{1}{1+x}$ \\
\hline
$r \leftarrow ($part of $R$ whose $e$-order is  $1) = -\xi$  \\
$[q_1', r'] \leftarrow$ $K[[x]]$-division($r, \inn_e(P_1), <_r$) \\
($q_1' = -\frac{1}{1+x}, r' = 0$) \\
$Q_1' \leftarrow -\frac{1}{1+x}$ \\
$R \leftarrow R - Q_1' P_1 = \frac{-1+x}{1+x}$ \\
\hline
$r \leftarrow ($part of $R$ whose $e$-order is $0) = \frac{-1+x}{1+x}$  \\
$[q_1', r'] \leftarrow$ $K[[x]]$-division($r, \inn_e(P_1), <_r$) \\
($q_1' = 0, r' = \frac{-1+x}{1+x}$) \\
$Q_1' \leftarrow 0$ \\
$R \leftarrow R - Q_1' P_1 = \frac{-1+x}{1+x}$ \\
\hline
\end{tabular}
\end{center}

The output of $\Dformal$-division($P, P_1$) is 
$$ \p^2 = (\frac{1}{1+x}\p - \frac{1}{1+x})((1+x)\p+x) +
   \frac{-1+x}{1+x}.$$ 
The quotient is  
$$ \frac{1}{1+x} \p - \frac{1}{1+x} 
   = (1-x+x^2-x^3+x^4-\cdots)(\p-1)$$
and the remainder is 
$$ \frac{-1+x}{1+x} = -1 + 2x - 2x^2 + 2x^3 - 2x^4 + 2x^5 - \cdots.$$

Next, we show the procedure of $\Dformal$-approximate-division
($P, P_1, 5$).

\begin{center}
\begin{tabular}{l}
\hline
$\overline{R} \leftarrow P, m_0 \leftarrow \ord_e(R) = 2$ \\
$M_0 \leftarrow 0, M_1 \leftarrow 1, M_2 \leftarrow 3$ \\
Bound $\leftarrow 5 + (0 + 1 + 3) = 9$ \\
\hline
$\overline{r} \leftarrow ($ part of $\overline{R}$ whose $e$-order is $2$ and 
total degree is less than $9$) $= \xi^2$  \\
$[\overline{q_1'}, \overline{r'}] \leftarrow$ $K[[x]]$-approximate-division(
$\overline{r}, \inn_e(P_1), <_r, 9$) \\
$(\overline{q_1'} = (1-x+x^2- \cdots + x^6) \xi, \overline{r'} = -x^7 \xi^2)$\\
$\overline{Q_1'} \leftarrow (1-x+x^2- \cdots + x^6) \p$ \\
$\overline{R} \leftarrow \overline{R} - \overline{Q_1'} P_1 = 
-x^7\p^2-\p-x^7\p-1+x-x^2+x^3-x^4+x^5-x^6$ \\
Bound $\leftarrow 9 - 3 = 6$ \\
\hline
$\overline{r} \leftarrow ($ part of $\overline{R}$ whose $e$-order is $1$ and 
total degree is less than 6$) = -\xi$  \\
$[\overline{q_1'}, \overline{r'}] \leftarrow$ $K[[x]]$-approximate-division(
  $\overline{r}, \inn_e(P_1), <_r, 6$) \\
$(\overline{q_1'} = -1+x-x^2+x^3-x^4, \overline{r'} = x^5 \xi)$\\
$\overline{Q_1'} \leftarrow -1+x-x^2+x^3-x^4$ \\
$\overline{R} \leftarrow \overline{R} - \overline{Q_1'} P_1 = 
  -x^7\p^2+x^5\p-x^7\p-1+2x-2x^2+2x^3-2x^4+2x^5-x^6$ \\
Bound $\leftarrow 6 - 1 = 5$ \\
\hline
$\overline{r} \leftarrow ($ part of $\overline{R}$ whose $e$-order is $0$ 
and total degree is less than $5)$\\ 
~~~~~~~ $= -1 + 2x - 2x^2 + 2x^3 - 2x^4$  \\
$[\overline{q_1'}, \overline{r'}] \leftarrow$ $K[[x]]$-approximate-division(
  $\overline{r}, \inn_e(P_1), <_r, 5$) \\
$(\overline{q_1'} = 0, \overline{r'} = \overline{r} = -1 + 2x - 2x^2 + 2x^3 - 2x^4)$ \\
$\overline{Q_1'} \leftarrow 0$ \\
$\overline{R} \leftarrow \overline{R} = 
-x^7\p^2+x^5\p-x^7\p-1+2x-2x^2+2x^3-2x^4+2x^5-x^6$\\
\hline
\end{tabular}
\end{center}

The output of $\Dformal$-approximate-division($P, P_1, 5$) is 
\begin{align*} 
\p^2 = & ((1 - x + x^2 + \cdots + x^6) \p + 
   (-1 + x - x^2 + x^3 - x^4)) \cdot ((1+x)\p+x) + \\
  & -x^7 \p^2 - x^7 \p + x^5 \p - 1 + 2x - 2x^2 + 2x^3 - 2x^4 +2x^5 - x^6
\end{align*}
The correct part of the remainder is 
$$ -1 + 2x - 2x^2 + 2x^3 - 2x^4.$$
And this part is the same with the part of the remainder of $\Dformal$-division
\end{eg}

\section{Computation of the Local $b$ Function by the Approximate Division Algorithm in $\Dformal[s]$}

We will apply the approximate division algorithm in $\Dformal$ (Algorithm \ref{ddivalgo}) 
to obtain the local $b$ function of a polynomial.
In the sequel, we use the parameter variable $s$ instead of $y$, and 
$\K$ denotes the field of complex numbers.
 
\subsection{Definition and Algorithm of the $b$ Function}
\begin{defn}
(Global $b$ function)

For given $f \in \K[x]$, 
we define the global $b$ function as the monic polynomial $\tilde{b}(s)$ of the least degree 
which satisfies $ \exists P \in D[s] \st P \cdot f^{s+1} = \tilde{b}(s) f^s$.
\end{defn}

\begin{defn}
(Local $b$ function)

For given $f \in \K[x]$,
we define the local $b$ function at the origin as the monic polynomial $b(s)$ of the least degree 
which satisfies $ \exists P \in \Dformal[s] \st P \cdot f^{s+1} = b(s) f^s$.
\end{defn}

\begin{eg}
($b$ function)

Let $f = x (x + y + 1) = x^2 + xy + x$.
The global $b$ function of $f$ is $(s+1)^2$, and $(\p_x^2 + \p_x \p_y) \cdot f^{s+1} = (s+1)^2 f^s$ holds.
The local $b$ function of $f$ is $s+1$, and $\frac{1}{1+2x+y} \p_x \cdot f^{s+1} = (s+1) f^s$ holds. 
\end{eg}

Oaku gave algorithms computing the global $b$ function and the local $b$ function
for a given polynomial 
(\cite{Oakubfunc1}, \cite{Oakubfunc2}, \cite{OTbfunc}).
A Gr\"{o}bner basis method in $D$ is used in these algorithm.

Noro gave an efficient algorithm computing the global $b$ function(\cite{Norobfunc}).
In this algorithm, a Gr\"{o}bner basis method in $D$, and modular method are used.
Especially, to eliminate variables, normal forms with respect to a Gr\"{o}bner basis are used.

In this paper, we propose an algorithm computing the local $b$ function by utilizing 
an approximate normal form with respect to a Gr\"{o}bner basis.
Our algorithm is an analogous algorithm with Noro's algorithm.
For this purpose, we review a Gr\"{o}bner basis and a normal form in $\Dformal$ 
and define an approximate normal form. 

\subsection{Gr\"{o}bner Basis and Normal Form in $\Dformal[s]$}

Although $\Dformal[s]$ is a transcedental object, 
for ideals in $\Dformal[s]$ generated by elements in $D[s]$, 
we can compute a Gr\"{o}bner basis in finite steps by either of the following algorithms.
\begin{itemize}
\item Lazard's method in $D$ (using a homogenization) 
\item method using the Mora division algorithm in $D$ (\cite{GO}, \cite{GOT})
\end{itemize}

\begin{defn}
(Normal form in $\Dformal[s]$)

For $P \in \Dformal[s]$ and a Gr\"{o}bner basis $G$ in $\Dformal[s]$,
the remainder of $\Dformal$-division($P, G$) (Theorem \ref{D-division}) is uniquely determined.
We call the remainder the normal form of $P$ by $G$ with respect to the monomial order $<$,
and we denote it by $\nf(P, G, <)$. 
\end{defn}

We note that a normal form does not have a finite representation in general.

\begin{lem}
(Ideal membership)

Let $G$ be a Gr\"{o}bner basis of an ideal $\mathcal{I}$ in $\Dformal[s]$. 
For $P \in \Dformal[s]$,
$P \in \mathcal{I}$ is equivalent to $\mathrm{NF}(P, G, <) = 0$.
\end{lem}

\begin{lem}
(Sum of normal forms)

Let $G$ be a Gr\"{o}bner basis in $\Dformal[s]$.
For $P, Q \in \Dformal[s]$, it holds that 
$$ \nf(P + Q, G, <) = \nf(P, G, <) + \nf(Q, G, <)$$
\end{lem}

\subsection{Algorithm Computing $\mathcal{I} \cap \K[s]$}
For an ideal $\mathcal{I}$ in $\Dformal[s]$, we propose an algorithm computing 
a generator of $\mathcal{I} \cap \K[s]$.
In this algorithm, we use approximate normal forms.
We fix a Gr\"{o}bner basis $G$ of $\mathcal{I}$.

For $g(s) = a_l s^l + a_{l-1} s^{l-1} + \cdots + a_0 \in \mathcal{I}$,
it holds that
\begin{align*}
g(s) \in \mathcal{I} &\Leftrightarrow \nf(g(s), G, <) = 0 \\
   &\Leftrightarrow 
   a_l \nf(s^l, G, <) + a_{l-1}\nf(s^{l-1}, G, <) + \cdots + a_0 \nf(1, G, <) = 0.
\end{align*}
We compute $\nf(s^i, G, <)$ beforehand, and
we compute the minimum $l$ and coefficients $a_l, \cdots, a_0 \in \K, a_l \neq 0$ which
satisfies
$$ a_l \nf(s^l, G, <) + a_{l-1}\nf(s^{l-1}, G, <) + \cdots + a_0 \nf(1, G, <) = 0.$$
Then the polynomial $a_l s^l + \cdots + a_0$ is a generator of the intersection.

However we cannot generally compute normal forms in $\Dformal[s]$ in finite steps.
The normal form $\nf(s^i, G, <)$ has infinite terms in general.
Therefore we use an approximate normal form.

\begin{defn}
(Approximate normal form)

Let $P \in \Dformal[s]$ and $G$ be a Gr\"{o}bner basis of an ideal $\mathcal{I}$ in $\Dformal[s]$.
We define the approximate normal form of $P$ by $G$ up to total degree $N$
as the remainder of $\Dformal$-approximate-division($P, G, N$) ~~(Algorithm \ref{ddivalgo}).
We denote it by $\nf(P, G, <, N)$.
By definition, 
$\nf(P, G, <)$ and $\nf(P, G, <, N)$ are the same up to the total degree $N - 1$.
\end{defn}

We suppose that 
$$ a_l \nf(s^l, G, <, N) + a_{l-1} \nf(s^{l-1}, G, <, N) + 
   \cdots + a_0 \nf(1, G, <, N) = 0.$$
Since we use approximate normal forms,
it is not always true that 
$a_l s^l + a_{l-1} s^{l-1} + \cdots + a_0 \in \mathcal{I}$.
However we can solve an ideal membership problem for algebraic data in $\Dformal[s]$
by utilizing the Mora division algorithm in $D[s]$.
Therefore we need to apply the Mora division algorithm to the gotten candidate 
$a_l s^l + \cdots + a_0$ in order to check if it is a member of $\mathcal{I}$.

We will explain an algorithm computing $\mathcal{I} \cap \K[s]$ by using these ideas.
We suppose $\mathcal{I} \cap \K[s] \neq \{0\}$.
We put 
\begin{align*}
 L_{i,N} &= \left\{(a_0, \cdots, a_i) \in \K^{i+1} ~|~ 
   a_i \nf(s^i, G, <, N) + \cdots + a_0 \nf(1, G, <, N) = 0 \right\} \\
 L_i     &= \left\{(a_0, \cdots, a_i) \in \K^{i+1} ~|~ 
   a_i \nf(s^i, G, <) + \cdots + a_0 \nf(1, G, <) = 0 \right\},
\end{align*}
where $G$ is a Gr\"{o}bner basis of $\mathcal{I}$ and $N \in \N$.
It holds that 
$$ L_{i,0} \supset L_{i,1} \supset L_{i,2} \supset \cdots \supset L_i.$$

The polynomial $a_i s^i + \cdots + a_0$ 
which satisfies $L_i \neq \{0\}, L_{i-1} = \{0\}$ 
and $(a_0, \cdots, a_i) \in L_i$ is 
a generator of $\mathcal{I} \cap \K[s]$.
In summary, an algorithm to compute the generator is as follows:

\begin{algo}
(Computing the generator of $\mathcal{I} \cap \K[s]$)
\label{find-gen}

\begin{enumerate}
\item[(1)] N $\leftarrow$ (a natural number), $d$ $\leftarrow$ 0
\item[(2)] Find $i$ which satisfies $L_{i-1,N} = \{0\}, L_{i,N} \neq \{0\}$.
      (Compute $L_{d,N}, L_{d+1,N}, \cdots $.
      Note that $L_{i-1,N} = \{0\} \Rightarrow L_{0} = L_{1} = \cdots = L_{i-1} = \{0\}$.)
\item [(3)] Take an element $(a_i, \cdots, a_0) \in L_{i,N}$ and 
      put $g(s) = a_i s^i + \cdots + a_0$.
\item [(4)] Divide $g(s)$ by $G$ by using the Mora division algorithm.
      If the remainder is $0$, then $g(s)$ is the generator.
      If not, then set $N \leftarrow N + 1, d \leftarrow i$ and goto (2).
\end{enumerate}
\end{algo}

\subsection{Algorithm Computing the Local $b$ Function}
For a polynomial $f \in \K[x]$,
an algorithm of computing the local $b$ function of $f$ at the origin 
is as follows:
\begin{enumerate}
\item[(1)] Compute a set of generators of $\mathrm{Ann}_{\widehat{\mathcal{D}}[s]} f^s$.
           Denote it by $H$.
\item[(2)] Set $\mathcal{I} = \Dformal[s] (H \cup \{f\})$.
      Compute the generator $b(s)$ of $\mathcal{I} \cap \K[s]$.
      The polynomial $b(s)$ is the local $b$ function.
\end{enumerate}

In (1), we take a set of generators of $\Ann_{D[s]} f^s$ as a set of generators of 
$\Ann_{\Dformal[s]} f^s$.
We can compute a set of generators of $\Ann_{D[s]} f^s$ by Oaku's algorithm
(\cite{Oakubfunc2}).
In (2), 
since $\mathcal{I} \cap \K[s] \neq \{0\}$ holds (from the existence of the local $b$ function), 
we can utilize the Algorithm \ref{find-gen}.

We show an example of computing the local $b$ function.
\begin{eg}
(Local $b$ function of $f=x^2(y+1)^2z^2$ at the origin)

The global $b$ function of $f$ is $(s+1)^3(s+\frac{1}{2})^3$, and 
the local $b$ function is $(s+1)^2(s+\frac{1}{2})^2$.

A set of generators $H$ of $\Ann_{\Dformal[s]} f^s$ is 
$$ \left\{P_1 = -2s + z \dz, 
     P_2 = - x \dx + z \dz, 
     P_3 = -\dy + x \dx  - y \dy \right\}$$
A Gr\"{o}bner basis $G$ of $\mathcal{J} = \Dformal[s] (H \cup \{f\})$ with respect to $<$ is 
$$
\left\{f, P_1, P_2, P_3, 
 P_4 = -x z^3 \dz - 2 x z^2, P_5 = -z^4 \dz^2 - 2 x z^2 \dx - 2 z^3 \dz - 2 z^2
\right\}  
$$
We set $N = 1$ and will apply the Algorithm \ref{find-gen} to compute the intersection $\mathcal{J} \cap \K[s]$.
We note that approximate normal forms of $s^i$ by $G$ up to total degree $N = 7$ are 
\begin{align*}
&\mathrm{NF}(1, G, <, 7) = 1 \\
&\mathrm{NF}(s, G, <, 7) = \frac{1}{2} z \dz \\
&\mathrm{NF}(s^2, G, <, 7) = \frac{1}{4} z^2 \dz^2 + \frac{1}{4} z \dz \\
&\mathrm{NF}(s^3, G, <, 7) = \frac{1}{8} z^3 \dz^3 + \frac{3}{8} z^2 \dz^2 
 + \frac{1}{8} z \dz\\
&\mathrm{NF}(s^4, G, <, 7) = -\frac{3}{8} z^3 \dz^3 - \frac{31}{16} z^2 \dz^2 
 - \frac{31}{16} z \dz - \frac{1}{4}.
\end{align*}
Steps of Algorithm \ref{find-gen} are 
\begin{center}
\begin{tabular}{ll}
$L_{1,1} = \{0\}$      &$L_1 = \{0\}$   \\
$L_{2,1} \ni (0,1)$    &$s \xrightarrow[G]{\text{Mora division}} \text{non-zero}$ \\
$L_{2,2} \ni (0,1)$    &$s \xrightarrow[G]{\text{Mora division}} \text{non-zero}$ \\
$L_{2,3} = \{0\}$      &$L_2 = \{0\}$   \\
$L_{3,3} \ni (0,-\frac{1}{2},1)$ &$s^2-\frac{1}{2}s \xrightarrow[G]{\text{Mora division}} \text{non-zero}$\\
$L_{3,4} \ni (0,-\frac{1}{2},1)$ &$s^2-\frac{1}{2}s \xrightarrow[G]{\text{Mora division}} \text{non-zero}$\\
$L_{3,5} = \{0\}$      &$L_3 = \{0\}$   \\
$L_{4,5} \ni (0, \frac{1}{2}, -\frac{3}{2}, 1)$ 
   &$s^3-\frac{3}{2}s^2+\frac{1}{2}s \xrightarrow[G]{\text{Mora division}} \text{non-zero}$  \\
$L_{4,6} \ni (0, \frac{1}{2}, -\frac{3}{2}, 1)$ 
   &$s^3-\frac{3}{2}s^2+\frac{1}{2}s \xrightarrow[G]{\text{Mora division}} \text{non-zero}$   \\
$L_{4,7} = \{0\}$      &$L_4 = \{0\}$   \\
$L_{5,7} \ni (\frac{1}{4}, \frac{3}{2}, \frac{13}{4}, 3, 1)$ 
   & $s^4+3s^3+\frac{13}{4}s^2+\frac{3}{2}s+\frac{1}{4} \xrightarrow[G]{\text{Mora division}} 0$.    
\end{tabular}
\end{center}
Therefore, the local $b$ function is
$s^4+3s^3+\frac{13}{4}s^2+\frac{3}{2}s+\frac{1}{4} = (s+1)^2(s+\frac{1}{2})$.
\end{eg}

\section{Implementation and Timing Data}
Our algorithm computing the local $b$ function has been implemented by utilizing computer algebra system 
``Risa/Asir''(\cite{RisaAsir}).
The name of our package is ``{\tt nk\_mora/local-b.rr}''.
We show the timing data.

Our algorithm computing the local $b$ function consists of the following 3 parts.
\begin{enumerate}
\item[(1)] Computation of a set of generators of $\Ann_{\Dformal[s]} f^s$ (In fact, computation of a set of generators of 
 $\Ann_{D[s]}f^s)$
\item[(2)] Computation of the Gr\"{o}bner basis of $\mathcal{I} = \Ann_{\Dformal[s]}f^s + \Dformal[s]f$
\item[(3)] Computation of the intersection $\mathcal{I} \cap \K[s]$ (Algorithm \ref{find-gen})
\end{enumerate}
To perform the step (1), we use Oaku's algorithm computing $\Ann_{D[s]}f^s$(\cite{Oakubfunc1}).
In the timing table, we denote this by $(\mathcal{I})$.
To perform the step (2), we have the following 2 algorithms.
\begin{enumerate}
\item[(2a)] Buchberger algorithm utlizing the Mora division algorithm in $D[s]$ 
\item[(2b)] Lazard's homogenized method in $D[s]$
\end{enumerate}
We denote this parts by (GB-a) and (GB-b).
To perform the step (3), we use the Algorithm \ref{find-gen}.
We denote this part by (localb-nf).

We took some examples from \cite{Oakubfunc2}.

{\scriptsize
\begin{center}
\begin{tabular}{|c|c|c|c|c|c|c|c|}
 \hline
 $f$         & $\mathcal{I}$& GB-a & GB-b    & localb-nf & localb & deg  
 \\
 \hline \hline
 $x+y^2+z^2$ & 0.00504   & 0.0190  & 0.00494 & 0.00503 & 0.00612 & 1 \\
 $x^2+y^2+z^2$ & 0.00524 & 0.00157 & 0.00224 & 0.0229  & 0.0119 & 2 \\
 $x^3+y^2+z^2$ & 0.00637 & 0.0298  & 0.0694 + 0.0192 & 0.0417+0.0244 & 0.0226 & 3 \\
 $x^4+y^2+z^2$ & 0.00631 & 0.0317  & 1.41 + 0.27 & 0.131 + 0.0884 & 0.0384 & 4 \\
 $x^5+y^2+z^2$ & 0.00624 & 0.0306  & 7.23 + 1.51 & 0.325 + 0.272  & 0.0612 & 5 \\
 $x^6+y^2+z^2$ & 0.00619 & 0.0291  & 33.0 + 3.25 & 0.726 + 0.598  & 0.129 & 6 \\
 $x^7+y^2+z^2$ & 0.00610 & 0.0281 + 0.0170 & 137.5 + 24.1 & 1.26 + 0.674 & 0.247 & 7 \\
 \hline
 $x^3+xy^2+z^2$ & 0.0260 & 0.172   & 0.107 & 0.114 + 0.0756  & 0.0767 & 4 \\  
 $x^4+xy^2+z^2$ & 0.0566 & 0.140   & 2.29 + 0.338 & 0.501 + 0.284 & 0.184 & 6 \\
 $x^5+xy^2+z^2$ & 0.0903 & 0.194   & 32.6 + 5.47  & 0.505 + 0.284 & 0.173 & 6 \\  
 \hline 
 $x^3+y^4+z^2$ & 0.00858 & 0.118   & 1.45 & 1.43 + 0.821 & 0.197 & 7 \\ 
 $x^3+xy^3+z^2$ & 0.0110 & 0.236 + 0.337 & 1.09 & 2.94 + 1.40 & 1.06 & 8 \\ 
 $x^3+y^5+z^2$ & 0.00979 & 0.123   & 18.3 + 2.82 & 4.36 + 1.80& 0.615 & 9\\
 \hline
 $x^5+x^3y^3+y^5$ & 0.70   & --- & 0.052 & --- & 7.06 & 7\\
 $x^4+x^2y^2+y^4$ & 0.0061 & 0.0138 & 0.0019 & 0.163 & 0.153 & 6 \\
 $x^3+x^2y^2+y^3$ & 0.031  & 0.0748 & 0.018  & ---   & 0.0932 & 4 \\
 $x^3+y^3+z^3$ & 0.0062    & 0.00162& 0.0025 & 0.142 & 0.112 & 5\\
 $x^6+y^4+z^3$ & 0.0134    & 0.052  & --- & --- & 23.2 & 18\\
 $x^3+y^2z^2$ & 0.029      & 0.0462 & 0.016  & 0.868 + 0.108 & 0.144 & 7\\
 $(x^3-y^2z^2)^2$ & 0.102  & 7.49 + 1.67 & 0.0932 & 17.4 + 1.40 & 3.99 & 14\\
 $x^5-y^2z^2$ & 0.013      & 0.0827 & 0.092  & 3.77 + 0.279 & 0.817 & 13 \\
 $x^5-y^2z^3$ & 0.0088     & 0.0265 & 0.0028 & 15.9 + 1.22 & 3.53 & 17\\
 $x^3+y^3-3xyz$ & 0.019    & 0.473  & 0.020  & 0.641 + 0.0540 & 0.0946 & 5\\
 $x^3+y^3+z^3-3xyz$ & 0.013& 0.140  & 0.012  & 0.0658 & 0.0207 & 3\\
 $y(x^5-y^2z^2)$ & 0.014   & 16.7 + 3.43 & 0.58 & 58.6 + 5.24 & 6.10 & 18\\
 $y(x^3-y^2z^2)$ & 0.644   & 2.65 + 0.668& 0.15 & 4.18 + 0.338 & 0.711 & 10\\
 $y((y+1)x^3-y^2z^2)$ & 0.284 & 18.8 + 1.34 & 3.6 + 0.14 & --- & 10.2 & 10\\   
 \hline
\end{tabular}
\end{center}
}
\begin{itemize}
\item localb  --- total time of Oaku's algorithm computing the local $b$ function (\cite{Oakubfunc1})
\item deg     --- degree of the local $b$ function
\item machine --- CPU : Athlon MP 1800+ (1533MHz) (2 CPU), Memory : 3GB, OS : FreeBSD 4.8
\end{itemize}

In the anterior half, Buchberger algorithm with Mora division (GB-a) is faster than Lazard homogenize method (GB-b).
While, in the last half, (GB-b) is faster than (GB-a).
From our experiments, in general, it seems that (GB-b) is faster than (GB-a).

In almost every case, Oaku's algorithm(localb) is faster than our algorithm($\mathcal{I}$ and GB-* and localb-nf).
The reason is that the computation of Gr\"{o}bner basis and division in $\Dformal$ is heavy.
Especially, the computation of $\Dformal$-division is heavy.

\end{document}